\documentclass[12pt]{article}
\usepackage{amsmath}
\usepackage{amssymb}
\usepackage{epic}
\usepackage{eepic}
\begin{document}

\newfont{\matha}{msbm10 scaled \magstep1}
\newfont{\mathb}{msbm8}
\newfont{\mathc}{eufm10 scaled \magstep1}

\newcommand{\Cbold}{\mbox{{\matha C}}}
\newcommand{\Nbold}{\mbox{{\matha N}}}
\newcommand{\Qbold}{\mbox{{\matha Q}}}
\newcommand{\Rbold}{\mbox{{\matha R}}}
\newcommand{\Zbold}{\mbox{{\matha Z}}}
\newcommand{\smCbold}{\mbox{{\mathb C}}}
\newcommand{\gbold}{\mbox{{\mathc g}}}
\newcommand{\nbold}{\mbox{{\mathc n}}}

\title{Auslander-Reiten Quivers and the Coxeter complex}
\author{Shmuel \textsc{ZELIKSON}
}
\date{August 7, 2002}
\maketitle

\begin{center}
\begin{minipage}{0.9\linewidth}
{\small \textsc{Abstract :} 
Let $Q$ be a quiver of type $ADE$. We construct the corresponding 
Auslander-Reiten quiver as a topological complex inside
the Coxeter complex 
associated with the underlying Dynkin 
diagram. We use the notion of chamber weights coming
from the theory of the canonical basis of
quantized enveloping algebras, and show this set has
a special linearity property in our setting. Finally, we consider
$A_n$ case, and describe how Auslander-Reiten 
quivers correspond to particular wiring diagrams.}
\end{minipage}
\end{center}

\section*{Introduction}

\hspace*{5mm} Let $D$ be a Dynkin diagram of type $ADE$, and $Q$ a quiver
obtained by orienting its edges. Associated with Q is a finite dimensional 
algebra over $\Cbold$, the path algebra $\Cbold Q$.
The category of representations of the quiver $Q$ over $\Cbold$ 
identifies with the category $\mbox{mod} \; \Cbold Q$ of finitely 
generated $\Cbold Q$ modules. Gabriel's theorem states
that  isomorphism classes of indecomposables ("indecomposables") of this category
are indexed by the positive roots of the root system
$\Phi$ associated with $D$. Furthermore only 
quivers of $ADE$ types have finitely many indecomposables.

The Auslander-Reiten quiver $\Gamma_Q$ codifies the structure of the category 
$\mbox{mod } \Cbold Q$. 
Vertices are the indecomposables, arrows are irreducible morphisms
between them. Our aim is to construct $\Gamma_Q$ as a topological complex
$\tilde{\Gamma}_Q$,
while retaining the spirit of Gabriel's theorem. We use the Coxeter 
complex $\Sigma$, that is the simplicial complex determined by
hyperplanes perpendicular to the positive roots, in order to get 
a unified labelling of both vertices and arrows of $\tilde{\Gamma}_Q$.

The vertices of $\tilde{\Gamma}_Q$ are one dimensional rays 
inside $\Sigma$, that is weights lying in Weyl group orbits of
the fundamental weights. To be more precise, $Q$ defines 
a bilinear form $R$ on the Euclidean space $E$ generated by $\Phi$
which corresponds to the homological form at the level of
$\mbox{mod} \; \Cbold Q$. Every column of the matrix of $R$
in the basis of simple roots  defines a weight $-\rho_i$.
The vertices of $\tilde{\Gamma}_Q$ obtain then out of the 
$\rho_i$'s  by the same linear combinations which construct
$\Phi^+$ out of the simple roots $\alpha_1, \,  \ldots, \, \alpha_n$.

$Q$ defines a particular reduced expression $\tilde{w}_0$
of the longest word of $W$, up to \linebreak[4]
commutation of simple reflections,
referred to as an adapted reduced expression.
$\tilde{w}_0$ induces a total order $\mathcal{R}_{\tilde{w}_0}$ 
of $\Phi^+$, and the main result of \cite{bedard} is
that the pair $\tilde{w}_0, \, \mathcal{R}_{\tilde{w}_0}$ 
encodes the structure of $\Gamma_Q$ (Theorem 1.1) . Conversely, one can
recover $\Gamma_Q$ out of $\tilde{w}_0$ (Proposition 1.2). 

Passing from the simple roots
to fundamental weights, one associates with $\tilde{w}_0$ 
the set of chamber weights $\mathcal{C}_{\tilde{w}_0}$ (Definition 2.1). 
This set was introduced in \cite{BZschubert} in the study
of total positivity questions in real semi-simple groups. 
We show (Corollary 2.5) $\mathcal{C}_{\tilde{w}_0}$ is the set 
of vertices defined above using  $R$. Actually one has linearity
of the $1:1$ correspondence $\varphi_{\tilde{w}_0} \; : \; 
\mathcal{R}_{\tilde{w}_0} \longrightarrow \mathcal{C}_{\tilde{w}_0}$, 
(Theorem 2.4), which turns out to be a restriction to $\Phi^+$
of a linear application $-\varphi_R$ defined by $R$.
We conjecture linearity is a characterization of  
reduced expressions $\tilde{w}_0$ adapted to quivers.

Two dimensional faces of $\Sigma$ (planar cones) are uniquely determined
by the two rays to which they are adjacent. They may therefore be seen
as edges between vertices. This allows us to use the set of vertices
$\mathcal{C}_{\tilde{w}_0}$ in order to associate 
with any reduced expression $\tilde{w}_0$ a quiver $\Gamma_{\tilde{w}_0}$
(Definition 3.5).
This quiver depends on $\tilde{w}_0$ only up to commutation
of simple reflections. When $\tilde{w}_0$ is adapted to a
quiver $Q$, $\Gamma_{\tilde{w}_0}$ depends only on $Q$, so we
can note it by $\tilde{\Gamma}_Q$. We show this quiver is 
isomorphic to $\Gamma_Q$ (Theorem 3.7). 
We use a direct combinatorial approach, thus 
isomorphism is proved independently of linearity. Note
however Theorem 2.4 defines naturally the isomorphism involved. 
The $\tilde{\Gamma}_{Q}$
being special cases of the $\Gamma_{\tilde{w}_0}$, we call the
later abstract Auslander-Reiten quivers.

The last section describes as an example $A_n$ type. 
This case is characterized by the fact all fundamental representations
are minuscule. We can therefore index \linebreak[4]
one-dimensional rays in $\Sigma$
by Young columns, and arrows of $\Gamma_{\tilde{w}_0}$ become
simply couples of column tableaus characterized by inclusion
of indices (Lemma 4.1).  This fact, combined with the existence of 
action of a Coxeter element in $\tilde{\Gamma}_Q$ (Proposition 2.7) similar
to the well-known one on $\Gamma_Q$, allows us to 
give a combinatorial description (Proposition 4.3) of
$\tilde{\Gamma}_Q$. On the theoretical level,
use of Young columns makes the link between 
$\Gamma_{\tilde{w}_0}$  and the wiring diagram $\mathcal{WD}(\tilde{w}_0)$ 
of $\tilde{w}_0$, used in \cite{BFZ} as a key tool
in the study of Lusztig's parameterizations of the canonical
basis (consult the more recent
\cite{BZinventiones} for use of chamber weights).

We are grateful to Bernhard Keller for helpful comments. 
We also wish to thank the referee for his remarks, 
which led us to improve the first version of this work, dealing
only with $A_n$ case.  This work was partially supported
by EC TMR network ``Algebraic Lie Representations'', 
grant no. ERB FMRX-CT97-01000.

\section{ Auslander-Reiten quivers}

\hspace*{5mm} We recall in this section basic facts about Auslander-Reiten
quivers, and proceed to study their link with reduced 
words of the longest word of the Weyl group. 

Let us fix $\Cbold$ as the ground field $k$.  Consider a $ADE$ 
Dynkin diagram $D$, of rank $n$, and index its vertices as in 
\cite{kac}, page 53. For sake of simplicity we shall identify
the set of simple roots $\{ \alpha_1 , \ldots \alpha_n \}$ 
with $I=\{ 1, \ldots n \}$. 

We note by $\Phi$ the root system associated with $D$, 
$E$ the real vector space in which it lies, 
$( \,, \,)$ the scalar product  defined on $E$ by the 
Cartan matrix, and by $\alpha_i^{\vee}$ the coroots. 
Let $W$ be the Weyl group, 
$\{ s_1, \, \ldots, \, s_n \}$ the set of simple reflections, $l()$ the
length function, and $w_0$ the unique longest element of $W$. 
$w_0$ induces an involution on $I$ 
(paragraphs (XI) of \cite{bourbaki}, Planches I,IV,V,VI,VII),
 which we shall note by $i \mapsto i^{*}$.

We obtain a quiver $Q$ 
by orienting the edges of $D$. 
$Q$ is given therefore as a couple $(Q_0,Q_1)$, with $Q_0$ ($=I$) a 
set of vertices, and $Q_1$ a set arrows. In this work
$Q$ will always denote a fixed quiver whose underlying 
Dynkin diagram is of type $ADE$.
Let $\mod \Cbold Q$ denote
the category of finite type representations of $Q$ over 
the complex numbers \cite{curtisreiner} 77.1.  
An object in this category 
is given by $M=(\underline{V}, \; \underline{\varphi})$,
where with every $i \in I$ we associate a finite dimensional
vector space $V_i$, and with every arrow $i \longrightarrow j$
in $Q_1$ a linear mapping $\varphi_{i,j} : V_i \longrightarrow V_j$.
One obtains the simple objects by associating a one dimensional
vector space to a vertex $i$, and $0$ spaces to the others. They are therefore
in $1:1$ correspondence with the simple roots.

We shall note by $\mbox{Ind} \; Q$ the set of isomorphism
classes $\lbrack M \rbrack$ of indecomposable \linebreak[4]
representations 
in $\mbox{mod} \; \Cbold Q$. The dimension
vector of a representation $M$ is the linear combination
$\mbox{dim} \; M=\underset{i \in I}{\sum} a_i \, \alpha_i$
with $a_i=\mbox{dim} \; V_i$.

\textbf{Gabriel's Theorem } Let $Q$ be a quiver of type $ADE$. 
Then the mapping \linebreak[4]
$ \lbrack M \rbrack \mapsto \mbox{dim} \; M$ establishes 
a one to one correspondence between $\mbox{Ind} \; Q$  and 
$\Phi^+$. \linebreak[4]
Furthermore, the quivers of type
$ADE$ are the only ones for which Ind $Q$ is a finite set.

\medskip

The category $\mbox{mod} \; \Cbold Q$  is hereditary. As observed
by Ringel \cite{ringelspecies} (see also \cite{ringelPBW}), 
the homological form
$$ \langle M , \, N \rangle := 
\dim_{\smCbold} \, \mbox{Hom} (M,N)-\dim_{\smCbold} \, \mbox{Ext}^{1}(M,N)$$
depends only on the dimension vectors. If we note
by $\lbrack \beta \rbrack$
the isomorphism class of the indecomposable corresponding to
$\beta \in \Phi^+$, then 
$\langle \lbrack \alpha \rbrack ,\, \lbrack \beta \rbrack \rangle=
R( \alpha , \, \beta )$ where $R$ is the bilinear form on $E$ defined
by 
$$ r_{i,j}:=R(\alpha_i, \; \alpha_j)=
\left\{
\begin{array}{ll}
 1 & \mbox{if} \; i=j \\
-1 & \mbox{if} \; i \longrightarrow j \in Q_1 \\
0 & \mbox{otherwise}
\end{array}
\right.$$
We shall refer to $R$ as the Ringel form of $Q$

The structure of $\mbox{mod} \; \Cbold Q$ is codified
by the Auslander-Reiten quiver $\Gamma_Q$. Vertices
are $\mbox{Ind} \; Q$, arrows are irreducible morphisms,
namely morphisms that cannot be written as a non-trivial composition
of two morphisms (\cite{smalo} page 166). $\Gamma_Q$ is endowed with 
translation  $\tau$  (\cite{smalo} page 225), which stratifies  
it into levels $\Gamma^i_{Q}$ indexed by $I$.  Such a level 
is of the form  
$$(P_i(Q)=)\tau^{N_i-1} I_i(Q) , \tau^{N_i-2} I_i(Q), \; \ldots \; , 
   \tau I_i(Q), I_i(Q)$$
with $P_i(Q)$ the projective cover of the simple module corresponding to $i \in I$, 
and $I_i(Q)$ the injective envelope of the simple module corresponding to $i^{*} \in I$.

A detailed construction of $\Gamma_Q$
our case is given  by \cite{gabriel} section 6.5. 
With $D$ one associates an infinite quiver $\Zbold D$
\cite{gabriel} fig 13. page 49, whose vertices
are $I \times \Zbold$. Arrows can exist between vertices
only if their first coordinates are linked in $D$.
The vertices of $D$ are ordered in a natural way such that
for $i,j$ linked in $D$, the corresponding arrows \linebreak[4] 
\vfill \eject
in $\Zbold D$
are given by the configuration

\medskip

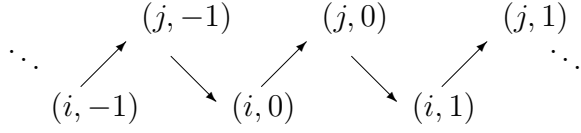
\unitlength=1mm
\begin{picture}(100,16)
\put(28,4){\vector(1,1){6}}
\put(40,8){\vector(1,-1){6}}
\put(52,4){\vector(1,1){6}}
\put(64,8){\vector(1,-1){6}}
\put(76,4){\vector(1,1){6}}
\put(18,6){$\ddots$}
\put(24,0){$(i,-1)$}
\put(36,12){$(j,-1)$}
\put(48,0){$(i,0)$}
\put(60,12){$(j,0)$}
\put(72,0){$(i,1)$}
\put(84,12){$(j,1)$}
\put(90,6){$\ddots$}
\end{picture}

\medskip

Thus $\Zbold D$ is an ``infinite mesh'' deduced out of $D$.
One defines on it the Nakashima involution $\nu$ 
\cite{gabriel} page 48.
The opposite quiver $Q^{op}$  to $Q$ may be realized inside $\Zbold D$ 
by a slice $S_{proj}$ starting at $(1,0)$. Its image by $\nu$
gives  a slice $S_{inj}$ which has the same arrows as $Q$.   
The Auslander-Reiten quiver then identifies with the subgraph 
of $\Zbold D$ delimited by $S_{proj}$ to the left, and $S_{inj}$ to
the right. Vertices of $S_{proj}$ (resp. $S_{inj}$) 
are the projectives in $\mbox{mod} \; \Cbold Q$ (resp. injectives). 
Each translation level $\Gamma^i_Q$ corresponds to the vertices in
its realization inside $\Zbold D$ having $i$ as first coordinate.
  One may deduce recursively the dimension vectors of the vertices
of $\Gamma_Q$ out of the slice $S_{proj}$
by use of additivity (\cite{gabriel} page 50)  
(for examples, see \cite{bedard} Examples 2.7, \; 2.8).

\textbf{Example 1.1 : }   
Consider the $A_n$ type quiver 
$Q^{(d)} \; : \; \stackrel{1}{\bullet} \longrightarrow  \stackrel{2}{\bullet} 
\longrightarrow \ldots \longrightarrow \stackrel{n}{\bullet}$, which
is distinguished by the regularity of its structure.
We reproduce the figure in \cite{ringelPBW} page 85, 
$\Gamma_{Q^{(d)}}$ is given by \hfill \break 
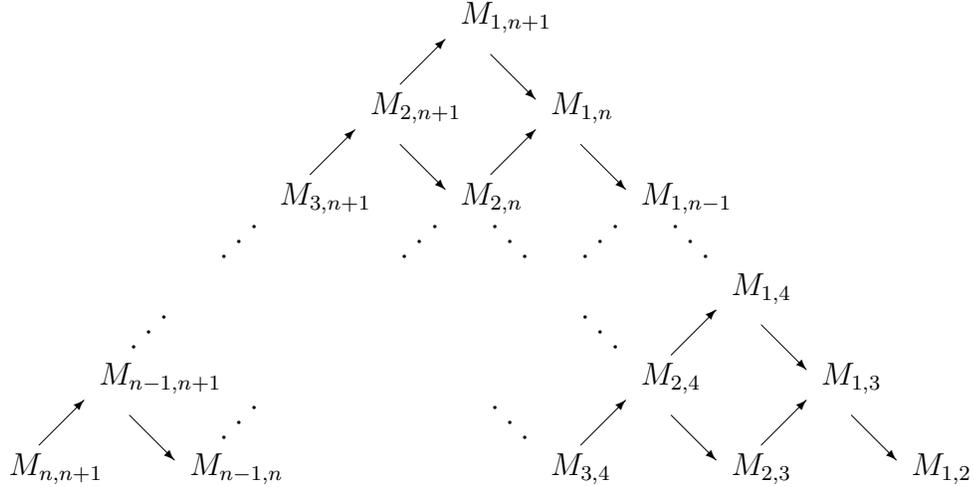
\medskip

\begin{picture}(140,60)(-15,0)
\put(4,4){\vector(1,1){6}}
\put(40,40){\vector(1,1){6}}
\put(52,52){\vector(1,1){6}}
\put(16,8){\vector(1,-1){6}}
\put(52,44){\vector(1,-1){6}}
\put(64,56){\vector(1,-1){6}}
\put(64,40){\vector(1,1){6}}
\put(76,44){\vector(1,-1){6}}
\put(76,4){\vector(1,1){6}}
\put(88,16){\vector(1,1){6}}
\put(88,8){\vector(1,-1){6}}
\put(100,20){\vector(1,-1){6}}
\put(100,4){\vector(1,1){6}}
\put(112,8){\vector(1,-1){6}}
\put(16,16){$\cdot$}
\put(18,18){$\cdot$}
\put(20,20){$\cdot$}
\put(28,28){$\cdot$}
\put(30,30){$\cdot$}
\put(32,32){$\cdot$}
\put(28,4){$\cdot$}
\put(30,6){$\cdot$}
\put(32,8){$\cdot$}
\put(52,28){$\cdot$}
\put(54,30){$\cdot$}
\put(56,32){$\cdot$}
\put(64,32){$\cdot$}
\put(66,30){$\cdot$}
\put(68,28){$\cdot$}
\put(76,28){$\cdot$}
\put(78,30){$\cdot$}
\put(80,32){$\cdot$}
\put(64,8){$\cdot$}
\put(66,6){$\cdot$}
\put(68,4){$\cdot$}
\put(76,20){$\cdot$}
\put(78,18){$\cdot$}
\put(80,16){$\cdot$}
\put(88,32){$\cdot$}
\put(90,30){$\cdot$}
\put(92,28){$\cdot$}
\put(0,0){$M_{n,n+1}$}
\put(12,12){$M_{n-1,n+1}$}
\put(36,36){$M_{3,n+1}$}
\put(48,48){$M_{2,n+1}$}
\put(60,60){$M_{1,n+1}$}
\put(24,0){$M_{n-1,n}$}
\put(60,36){$M_{2,n}$}
\put(72,48){$M_{1,n}$}
\put(84,36){$M_{1,n-1}$}
\put(72,0){$M_{3,4}$}
\put(84,12){$M_{2,4}$}
\put(96,24){$M_{1,4}$}
\put(96,0){$M_{2,3}$}
\put(108,12){$M_{1,3}$}
\put(120,0){$M_{1,2}$}
\end{picture}
\vspace*{5mm}
\hfill \break
where 
$M_{i,j+1}:=\lbrack \alpha_i + \alpha_{i+1}+ \ldots + \alpha_j \rbrack$.
We have 
$S_{proj}=\{ M_{1,n+1}, \, M_{2,n+1} \ldots M_{n,n+1}\}$,  
($M_{n,n+1}$ corresponding to $(0,1) \in \Zbold D$), 
and $S_{inj}=\{ M_{1,n+1}, \, M_{1,n} \ldots M_{1,3}, \, M_{1,2} \}$.
For $i \in  I$, the translation level 
$\Gamma^i_{Q^{(d)}}$ is $i^{\mbox{\small th}}$ horizontal level
from the top (i.e. containing $\lbrack M_{i,n+1} \rbrack$).

Given $w \in W$, we shall use the notation  
$\tilde{w}=s_{i_1} \ldots s_{i_m}$ for a reduced expression of $w$, 
the $\sim$ indicating we are referring not to $w$ but to its reduced 
expression. Recall one has for the longest element $l(w_0)=N$ where  
$N=\mbox{card} \Phi^+$. A fixed reduced 
expression $\tilde{w_0}=s_{i_1} s_{i_2} \ldots s_{i_N}$
induces a total ordering of the positive
roots $\mathcal{R}_{\tilde{w}_0}=(\beta_1, \beta_2 \ldots \beta_N)$ by
$$ \beta_j= s_{i_1} s_{i_2} \ldots s_{i_{j-1}} ( \alpha_{i_j} ).$$
 
$\mathcal{R}_{\tilde{w}_0}$ obtains as a linear refinement of 
convex partial order on $\Phi^+$. Convex (or normal) 
means a partial order $\preccurlyeq$
verifying conditions \cite{BFZ} :
$$ \begin{array}{cc}
\forall \alpha, \beta \in \Phi^+, 
\; \mbox{such that} \; \gamma=\alpha+\beta \in \Phi^+ \; : \\
\mbox{either} \; 
\alpha \preccurlyeq \gamma \preccurlyeq \beta
\; \mbox{or} \;  
\beta \preccurlyeq \gamma \preccurlyeq \alpha. 
\end{array}$$

Conversely, a convex order $\preccurlyeq$, 
determines  $\tilde{w}_0$, up to commutation of simple reflections, 
that is up to "2-moves" 
replacing $s_{i_j} s_{i_{j+1}}$ by $s_{i_{j+1}} s_{i_j}, 
\; i_j$  and $i_{j+1}$
not linked in $D$. If we consider two reduced expressions of $w_0$ 
as equivalent
if one may be obtained from the other by a sequence of 2-moves, 
then the equivalence classes we obtain are known as commutation  
classes \cite{bedard} 1.4. We have  a $1:1$ correspondence 
between them and convex orderings. All constructions in the 
sequel will actually depend only on the commutation class
(or equivalently a convex ordering). 

A quiver $Q$ defines a convex ordering $\preccurlyeq_Q$ \cite{ringelPBW}. One may use
the Ringel form, putting $\beta \preccurlyeq \gamma$ if 
$\langle \beta, \, \gamma \rangle > 0$, and taking the transitive 
closure of these relations. The same partial order may be obtained
from $\Gamma_Q$, one has $\alpha \preccurlyeq \beta$ iff
there is a path in $\Gamma_Q$ from  $\lbrack \alpha \rbrack$
to $\lbrack \beta \rbrack$.
Equivalence of the two constructions follows form the fact
arrows in $\Gamma_Q$ are 
couples $\lbrack \alpha \rbrack, \, \lbrack \beta \rbrack$ such that 
$ R (\alpha , \, \beta ) > 0$ and there is no $\gamma$ 
such that $R( \alpha , \, \gamma) > 0$ and
$R (\gamma, \, \beta) >0$.
We shall say $\preccurlyeq$ is the convex  ordering adapted to $Q$, 
and use the same reference to reduced expressions in the corresponding commutation
class. 

We say, following \cite{BGP} a vertex $k$ 
is a sink
of $Q$ if there are only entering arrows into
it. Call it a source if there are only exiting
arrows out of it. If $k$ is a sink of $Q$
define $S_k Q$ as the quiver obtained by
reversing the direction of the arrows with
$k$ as end, thus turning $k$ into a source.
A reduced expression 
$\tilde{w}_0=s_{i_1} s_{i_2} \ldots s_{i_N}$ 
is adapted to $Q$ if and only if $i_1$ is a sink of $Q$, $i_2$ a sink of 
$S_{i_1} Q$, $i_3$ is a sink of $S_{i_2} S_{i_1} Q$ and so on.

\textbf{Example 1.2 : }
Consider $Q=Q^{(d)} : \; \;  \stackrel{1}{\bullet}
\longrightarrow \stackrel{2}{\bullet} 
\longrightarrow \stackrel{3}{\bullet} $ of type $A_3$. 
$1$ is a source, $3$ is a sink, and $Q'=S_3 Q$ is then 
$$ Q' : \; \; \stackrel{1}{\bullet}
\longrightarrow \stackrel{2}{\bullet} 
\longleftarrow \stackrel{3}{\bullet} $$ 
 
We leave it to the reader to verify
$\tilde{w}_0=s_3 s_2 s_1 s_3 s_2 s_3$ is adapted
to $Q^{(d)}$ while $\tilde{w}_0=s_2 s_1 s_3 s_2 s_3 s_1$ is
adapted to $Q'$. For arbitrary $n$, one may verify \linebreak[4]
$\tilde{w}_0=s_n s_{n-1} \ldots s_1 s_n s_{n-1} \ldots s_2 
\; \ldots \; s_n s_{n-1} s_n$
 is adapted to the distinguished orientation $Q_d$.

Given $Q$, let us renumber its vertices in a compatible
way with its arrows. Namely we renumber $\{ \alpha_1 , \, 
\ldots \, \alpha_n \}$ as
$\{ \beta_1, \, \ldots , \,  \beta_n \}$ such that
$\beta_j \longrightarrow \beta_k$ in $Q$ implies
$j>k$. One then associates to Q (\cite{BGP} \S 1.2) a Coxeter
element $c=s_{\beta_n} \ldots s_{\beta_2} s_{\beta_1}$. 
There might be several compatible indexations, yet
by \cite{BGP} Lemma 1.2, $c$ is independent of
a particular choice. For instance 
$c=s_3 s_2 s_1$ is the Coxeter element associated 
to the quiver $Q^{(d)}$ in the example above, 
and $c=s_1 s_3 s_2 (=s_3 s_1 s_2)$ the Coxeter element
to the quiver $Q'$. 
Observe, as in the proof of \cite{BGP} Theorem 1.2, 
$s_{\beta_n} \ldots   s_{\beta_2} s_{\beta_1}$ is
a sequence adapted to $Q$.
It is well known  
that translation $\tau$ of $\Gamma_Q$ becomes
action of the Coxeter element $c$ at the level
of dimension vectors.

If $k$ is a sink of $Q$, and $Q':=S_k Q$
then the categories $\mbox{mod} \; \Cbold Q$ and $\mbox{mod} \; \Cbold Q'$
 are closely linked  by the BGP-reflection
functor \cite{curtisreiner} 77.4 corresponding to $s_k$.
Bedard \cite{bedard} combines use of these functors together 
with the description of $\Gamma_Q$ inside $\Zbold D$ in order to
show how to obtain from $\Gamma_Q$ all $\tilde{w}_0$
adapted to $Q$.

If for $i \in I$, the occurrences of $s_i$ in $\tilde{w}_0$ 
are at positions $j_1, \, j_2 \ldots j_r$  define 
$\mathcal{R}^i_{\tilde{w}_0}$ as the ordered $r$-tuple
$(\gamma_1,  \, \gamma_2, \ldots , \, \gamma_r)$
with $\gamma_k=s_{i_1} s_{i_2} \ldots s_{i_{j_k-1}} (\alpha_{j_k})$
Theorem 2.17 of \cite{bedard} may
\vfill \eject 

then be restated as follows : 

\medskip

\textbf{Theorem 1.1}  Let $Q$ be a quiver of type $ADE$.

\textbf{a)} Associate with every vertex $\lbrack \beta  \rbrack$ of  
$\Gamma_Q$ the translation level
$i_{\beta}$ it belongs to. Read $\mbox{Ind} \;  Q$ sequentially, 
in a manner compatible with the arrows 
(i.e $\lbrack \alpha \rbrack \longrightarrow \lbrack \beta \rbrack$ implies
$\lbrack \alpha \rbrack$ appears before $\lbrack \beta \rbrack$). Replacing
vertices $\lbrack \beta \rbrack$ by the corresponding $i_{\beta}$ one
produces an $N$-tuple $(i_1, i_2, \ldots i_N)$ giving
a reduced expression adapted to $Q$. Furthermore, all $\tilde{w}_0$
adapted to $Q$ are obtained this way.

\textbf{b)} Fix $\tilde{w}_0$ as in  a). For all $i \in I$ taking dimension
vectors establishes a  $1:1$
correspondence (as ordered sets) between $\Gamma_Q^i$ and 
$\mathcal{R}^i_{\tilde{w}_0}$. 
\medskip

Part $b)$ tells us in particular each $\lbrack \alpha \rbrack \in \Gamma^i_Q$
has the form $\alpha=s_{i_1} s_{i_2} \ldots s_{i_{j-1}} (\alpha_{i_j})$ 
with $i_j=i$. We shall note by $N_i(Q)$ the  cardinality 
of $\Gamma^i_Q$ (or of $\mathcal{R}^i_{\tilde{w}_0}$).

\textbf{Example 1.3 : }
Consider the quiver $Q:= \stackrel{1}{\bullet}
\longrightarrow \stackrel{2}{\bullet} 
\longleftarrow \stackrel{3}{\bullet}$ of type $A_3$. 
$\Gamma_Q$ is

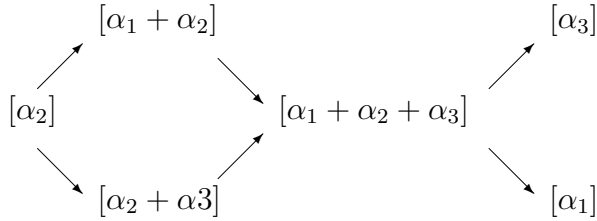
\vspace*{-1.9cm}
\begin{picture}(80,50)
\put(16,16){\vector(1,1){6}}
\put(40,20){\vector(1,-1){6}}
\put(76,16){\vector(1,1){6}}
\put(16,8){\vector(1,-1){6}}
\put(40,4){\vector(1,1){6}}
\put(76,8){\vector(1,-1){6}}
\put(12,12){$\lbrack \alpha_2 \rbrack$}
\put(24,24){$\lbrack \alpha_1+\alpha_2 \rbrack$}
\put(24,0){$\lbrack \alpha_2+\alpha3 \rbrack$}
\put(48,12){$\lbrack \alpha_1+\alpha_2+\alpha_3 \rbrack$}
\put(84,24){$\lbrack \alpha_3 \rbrack$}
\put(84,0){$\lbrack \alpha_1 \rbrack$}
\end{picture}
\bigskip

Replacing vertices by their translation levels one obtains 

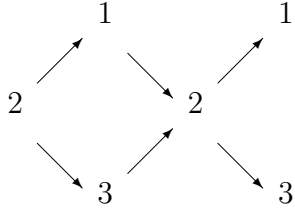
\vspace*{-1.9cm} \hspace*{2cm} 
\begin{picture}(100,50)
\put(16,16){\vector(1,1){6}}
\put(28,20){\vector(1,-1){6}}
\put(40,16){\vector(1,1){6}}
\put(16,8){\vector(1,-1){6}}
\put(28,4){\vector(1,1){6}}
\put(40,8){\vector(1,-1){6}}
\put(12,12){$2$}
\put(24,24){$1$}
\put(24,0){$3$}
\put(36,12){$2$}
\put(48,24){$1$}
\put(48,0){$3$}
\end{picture}
\bigskip

There are $4$ ways of reading the vertices in a compatible way
with the arrows  \linebreak[4]
$(2,1,3,2,1,3), \; (2,1,3,2,3,1), \; (2,3,1,2,1,3), \; (2,3,1,2,3,1)$.
The $4$ corresponding reduced expressions form the commutation class
determined by $Q$.

Theorem 1.1 allows us also to describe the opposite direction, recovering
$\Gamma_Q$ from a $\tilde{w}_0$ adapted to $Q$. Citing the referee, 
this result is well known to the experts, yet unpublished so we provide
its proof.

\medskip

\textbf{Proposition 1.2} Let $\tilde{w}_0$ be a reduced expression 
adapted to a quiver $Q$ of type $ADE$. Consider a couple of roots 
$\alpha=s_{i_1} s_{i_2} \ldots s_{i_{j-1}} (\alpha_{i_j}), \;
\beta=s_{i_1} s_{i_2} \ldots s_{i_{k-1}} (\alpha_{i_k})$ in 
$\mathcal{R}_{\tilde{w}_0}$. 
$ \lbrack \alpha \rbrack \longrightarrow \lbrack \beta \rbrack$ 
is an arrow in $\Gamma_Q$ iff 
\begin{itemize}
\item[i)] $i_j$ and $i_k$ are linked in the Dynkin diagram $D$,
\item[ii)] $j<k$,
\item[iii)] For $j < p < k$, $i_p \neq i_j$.
\end{itemize}

\vfill \eject
\textbf{Proof : } Since $\Gamma_Q$ lies
inside $\Zbold D$, observe there are arrows only between translation
levels whose indices are linked in $D$. We need therefore 
only to consider
the case of $i_j$ and $i_k$ linked in $D$. Put $a:=i_j, \; b:=i_k$.
As $\tilde{w}_0$ is adapted to $Q$, by Theorem 1.1 it
obtains by reading translation levels of vertices of $\Gamma_Q$
in a compatible manner with its arrows. Now the mesh structure
of $\Zbold D$ implies an arrow from a vertex 
$\lbrack \alpha \rbrack$ on  $\Gamma^a_Q$ to a vertex $\lbrack \beta \rbrack$
on $\Gamma^b_Q$ is immediately followed by an arrow from
$\lbrack \beta \rbrack$ to the vertex $\lbrack \gamma \rbrack$ following
$\lbrack \alpha \rbrack$ on the translation level $\Gamma^a_Q$, and vice-versa.
In order to comply with this configuration, if we delete from $\tilde{w}_0$
all simple reflection $s_{i_l}$ with $l \neq a,b$
we have to obtain a word either of the form
$s_a s_b s_a s_b \ldots$ or of the form $s_b s_a s_b s_a \ldots$. 
By  ii) and iii), our criterion gives then the same structure of arrows as the 
one of $\Gamma_Q$ inside $\Zbold D$ $\square$.

\medskip

\textbf{Definition 1.3 } Let $\tilde{w}=s_{i_1} \ldots s_{i_m}$ 
be a reduced expression
of $w \in W$. We shall say $\tilde{w}$ is alternating if
for every couple $j, \, k$ of indices linked in $D$, the word one obtains
by erasing from $\tilde{w}$ all $s_{i_r}, \; i_r \neq j,k$ 
is either  of the form $s_i s_j s_i \ldots$ or of the form
$s_j s_i s_j \ldots$. 

\medskip 

\textbf{Lemma 1.4 } Any reduced expression $\tilde{w}_0$ adapted to a quiver $Q$ is
alternating.

The proof is included in that of Proposition 1.2.

\textbf{Example 1.4 : } Consider $\tilde{w}_0=s_3 s_2 s_1 s_3 s_2 s_3$.
For the couple $\{ 1, \, 2 \}$, we get by erasing occurrences of 
$s_3$, the word $s_2 s_1 s_2$, and for the couple
$\{ 2 , \, 3 \}$ erasing occurrences of $s_1$ gives
$s_3 s_2 s_3 s_2 s_3$.

\section{The linearity phenomenon}

\hspace*{5mm} We import in this section the notion of chamber weights 
$\mathcal{C}_{\tilde{w}_0}$ from quantum groups  and show in case $\tilde{w}_0$
is adapted to $Q$, the correspondence 
$\varphi_{\tilde{w}_0} \, : \, \mathcal{R}_{\tilde{w}_0} 
\longrightarrow \mathcal{C}_{\tilde{w}_0}$ 
is linear and given by the Ringel form. 

Let us consider orbits of fundamental weights 
$\omega_i$ instead of positive roots.

\medskip
 
\textbf{Definition 2.1 }
Let $W$  be a Weyl group of finite type, and let us fix 
a reduced expression $\tilde{w}_0$.   
The set of chamber weights associated with $\tilde{w}_0$
is the $N$-tuple \linebreak[4]
$\mathcal{C}_{\tilde{w}_0}:=(\mu_1, \, \mu_2 , \, \ldots , \, \mu_N)$
where $\mu_j=s_{i_1} s_{i_2} \ldots s_{i_{j}} (\omega_{i_j})$.

 Chamber weights were introduced in \cite{BZschubert} as a key tool
in the study of total positivity in semi-simple groups. Our convention
is in opposite order to that in \cite{BZschubert}. 
All elements of $\mathcal{C}_{\tilde{w}_0}$ 
are distinct (\cite{BZschubert} 2.9).
By its construction, $\mathcal{C}_{\tilde{w}_0}$ seen as a set  
depends only on the commutation class of $\tilde{w}_0$.
Let us observe the cardinality of $\mathcal{F}=\underset{i \in I}{\bigcup} W \omega_i$
is greater than $N$, so unlike the positive roots case, $\tilde{w}_0$ actually
``cuts out'' a subset inside $\mathcal{F}$.

For any $\tilde{w}_0$, one has a natural $1:1$ correspondence 
$\varphi_{\tilde{w}_0}$ between 
$\mathcal{R}_{\tilde{w}_0}=(\beta_1, \,  \, \ldots , \, \beta_N)$
and $\mathcal{C}_{\tilde{w}_0}=(\mu_1, \, \ldots, \, \mu_N)$   defined by 
$\varphi_{\tilde{w}_0}(\beta_k)=\mu_k, \; k=1, \ldots, \, N$.
Observe this bijection depends only on the commutation class, as 
interchange of two
commuting reflections $s_{i_j}, \, s_{i_{j+1}}$ results in an interchange 
between $\beta_j, \, \beta_{j+1}$ in $\mathcal{R}_{\tilde{w}_0}$, and 
$\mu_j, \; \mu_{j+1}$ in $\mathcal{C}_{\tilde{w}_0}$.

If we consider columns of the matrix $(r_{i,j})$ of the
 Ringel form $R$ as coordinate vectors in the basis of the $\omega_i$'s, we
obtain the following weights :     
$$\forall i \in I \; : \; \rho_i:=\underset{l \in I}{\sum} r_{l,i} \omega_i.$$
We note by $\varphi_R$ the linear mapping from $E$ to $E$ defined by
$ \varphi_R(\alpha_i)=\rho_i$ for $i$ in $I$.

\medskip

\textbf{Lemma 2.2}  Suppose $k$ is a sink of $Q$ and  
$Q':=S_k Q$. Let $R$ and $R'$ be the corresponding Ringel
forms. Then 
$$ s_k \varphi_{R'}=\varphi_R s_k $$
\textbf{Proof : } The lemma is a direct consequence of the
behavior of Ringel forms under $BGP$ reflection functors. 
Let $\{ \rho_i \}_{i \in I}, \; \{ \rho'_i \}_{i \in I}$ be the
respective weights defined by $R$ and $R'$.
Note by $D_k$ the set of $j \neq k$ linked to $k$ in $D$.

The coefficient of $\omega_k$ in $\rho_j$ verifies 
$\lbrack \rho_j \, : \, \omega_k \rbrack=0$ for $j \neq k$, while
$\rho_k=\omega_k -\underset{j \in D_k}{\sum} \omega_j$.
In comparison, we have for $Q'$ : 
$$\rho'_j=
\begin{cases}
\omega_k & j=k, \\
\rho_j -\omega_k & j \in D_k, \\
\rho_j & \mbox{otherwise}.  
\end{cases}$$ 

Let us show directly 
$\forall j \in I \; : \; \varphi_{R'}(\alpha_j)= s_k \varphi_R s_k (\alpha_j)$
(we give details for convenience of the reader). 

Case 1 : $j=k$.  

$$ \begin{array}{rcl} 
s_k \varphi_R s_k (\alpha_k) & = & s_k(-\rho_k) \\
& = & s_k (-\omega_k+\underset{j \in D_k}{\sum} \omega_j) \\
&= & -w_k+\alpha_k+\underset{j \in D_k}{\sum} \omega_j.
\end{array}$$
Now $\alpha_k=2 \omega_k-\underset{j \in D_k}{\sum} \omega_j$ so
we get $s_k (-\rho_k)=\omega_k=\rho'_k=\varphi_{R'}(\alpha_k)$.

Case 2 : $j \in D_k$.
$$ \begin{array}{rl}
 s_k \varphi s_k (\alpha_j) & =s_k \varphi_R(\alpha_j+\alpha_k) \\
  & =s_k (\rho_j+\rho_k).
\end{array}$$
$s_k (\rho_j)=\rho_j$ since $\lbrack \rho_j \, : \, \omega_k \rbrack=0$. 
We have just seen $s_k (\rho_k)=-\omega_k$. We get therefore 
$s_k \varphi_R s_k (\alpha_j)=\rho_j-\omega_k=\rho'_j=\varphi_{R'}(\alpha_j)$

Case 3 : $j \neq k$ and not in $D_k$.

Immediate, $\rho_j=\rho'_j$ and $s_k$ stabilizes both $\alpha_j$ and $\rho_j$
 $\square$

\medskip

\textbf{Corollary 2.3} Assume the 
expression $s_{i_1} s_{i_2} \ldots s_{i_{k-1}}$ 
is adapted to $Q$. Note for 
$j=1 \ldots k-1 \; : \; Q_j:=S_{i_{j-1}} S_{i_{j-2}} \ldots S_{i_1} Q \; 
(Q_1:=Q)$
and $R_j$ the Ringel form corresponding to $Q_j$. Then
$$ \varphi_{R_1} s_{i_1} s_{i_2} \ldots s_{i_{k-1}}= 
    s_{i_1} s_{i_2} \ldots s_{i_{k-1}} \varphi_{R_k}.$$

\textbf{Proof : }
For any $j=1, \ldots k-1$, $i_j$ is a sink of $Q_j$ so
that $\varphi_{R_{j}} s_{i_j} =s_{i_j} \varphi_{R_{j+1}}$.

\medskip

\textbf{Theorem 2.4} (Linearity) Let $\tilde{w}_0$ 
a reduced expression adapted to $Q$ of type $ADE$.
Then the mapping $\varphi_{\tilde{w}_0}$ is linear and
given by the Ringel form $R$: 
$$ \forall  1 \leq k \leq N, \; \mu_k=-\varphi_R(\beta_k)$$

\medskip

\textbf{Corollary 2.5} The $-\rho_i$'s, 
$i \in I$ are elements of  $\mathcal{C}_{\tilde{w}_0}$, and
  $\mathcal{C}_{\tilde{w}_0}$ is obtained from them by the 
same linear combinations
which construct $\Phi^+$ out of the simple roots $\{ \alpha_1 , \ldots , \, \alpha_n \}$.

\medskip

\textbf{Proof of Theorem 2.4 : }
Consider an arbitrary $k, \; 1 \leq k \leq N$. By hypothesis, 
$\tilde{w}_0=s_{i_1} s_{i_2} \ldots s_{i_N}$ is adapted to $Q$ 
so may apply Corollary 2.3 and its conventions (in particular $R_1=R$).
$$ \begin{array}{rl} 
\varphi_{R} (\beta_k) &= \varphi_{R_1} s_{i_1} s_{i_2} \ldots s_{i_{k-1}} (\alpha_{i_k}) \\ 
&= s_{i_1} s_{i_2} \ldots s_{i_{k-1}} \varphi_{R_k} (\alpha_{i_k}) 
\\
&= s_{i_1} s_{i_2} \ldots s_{i_{k-1}} \varphi_{R_k} (s_{i_k} (-\alpha_{i_k}))
\\
&=s_{i_1} s_{i_2} \ldots s_{i_{k-1}} s_{i_k} \varphi_{R_{k+1}} (-\alpha_{i_k})
\\
&=-s_{i_1} s_{i_2} \ldots s_{i_{k-1}} s_{i_k} \varphi_{R_{k+1}} (\alpha_{i_k})
\end{array}$$

Observe $i_k$ is a source of $Q_{k+1}$. We have seen
while proving Lemma 2.2, in that case 
$\varphi_{R_{k+1}}(\alpha_{i_k})=\omega_{i_k}$,
so finally 
$$\varphi_R(\beta_k)=-s_{i_1} s_{i_2} \ldots s_{i_k} (\omega_{i_k})=-\mu_k. \; \; \; \square$$

\textbf{Examples 2.1 } Consider $A_3$ type.

\textbf{1)} We have seen $\tilde{w}_0=s_2 s_1 s_3 s_2 s_3 s_1$ is adapted
to a quiver. $\mathcal{R}_{\tilde{w}_0}$ is given by
$$ \begin{array}{lll}
\beta_1=\alpha_2, & \hspace*{4mm} \beta_3=\alpha_2+\alpha_3, 
& \hspace*{4mm} \beta_5=\alpha_1, \\
\beta_2=\alpha_1+\alpha_2, & \hspace*{4mm} \beta_4=\alpha_1+\alpha_2+\alpha_3, & 
\hspace*{4mm} \beta_6=\alpha_3, \end{array}$$
and $\mathcal{C}_{\tilde{w}_0}$ by
$$ \begin{array}{lll}
\mu_1=\omega_1-\omega_2+\omega_3, & \hspace*{4mm} \mu_3=\omega_1-\omega_2, & 
\hspace*{4mm} \mu_5=-\omega_1, \\
\mu_2=-\omega_2+\omega_3, & \hspace*{4mm} \mu_4=-\omega_2, & 
\hspace*{4mm} \mu_6=-\omega_3. 
\end{array}$$
The matrix of the Ringel form is 
$$ \left\lbrack
\begin{array}{rrr} 
1 & -1 & 0 \\
0 & 1 & 0 \\
0 & -1 & 1 
\end{array} 
\right\rbrack$$
$\mu_5, \, \mu_1, \, \mu_6$ are the respective images by $\varphi_{\tilde{w}_0}$
of the simple roots, and indeed \linebreak[4]
$-\rho_1=-\omega_1=\mu_5$, 
$-\rho_2=\omega_1-\omega_2+\omega_3=-\mu_1$, 
$-\rho_3=-\omega_3=\mu_6$.
Finally, observe one has $\mu_2=\varphi_{\tilde{w}_0} (\alpha_1+\alpha_2)$,
and as expected $-\rho_1-\rho_2=-\omega_2+\omega_3$. Similarly, 
$\mu_3=-\rho_2-\rho_3$ and $\mu_4=-\rho_1-\rho_2-\rho_3$.

\textbf{2)} The reduced expression $\tilde{w}_0=s_2 s_1 s_2 s_3 s_2 s_1$ 
is not adapted to a quiver.  \linebreak[4]
$\mathcal{R}_{\tilde{w}_0}=(\beta_1' , \, \ldots , \, \beta_6')$ with
$$ \begin{array}{lll}
\beta_1'=\alpha_2, & 
\hspace*{4mm} \beta_3'=\alpha_1, & \hspace*{4mm} \beta_5'=\alpha_2+\alpha_3, \\
\beta_2'=\alpha_1+\alpha_2, & \hspace*{4mm} \beta_4'=\alpha_1+\alpha_2+\alpha_3, &
\hspace*{4mm} \beta_6'=\alpha_3, 
\end{array}$$
while $\mathcal{C}_{\tilde{w}_0}=(\mu_1', \, \ldots, \, \mu_6')$ with
$$ \begin{array}{lll}
\mu_1'=\omega_1-\omega_2+\omega_3, & \hspace*{4mm} \mu_3'=-\omega_1+\omega_3, & 
\hspace*{4mm} \mu_5'=-\omega_2, \\
\mu_2'=-\omega_2+\omega_3, 
& \hspace*{4mm} \mu_4'=-\omega_1, & \hspace*{4mm} \mu_6'=-\omega_3.
\end{array}
$$
In this case, linearity fails. Observe
$\mu_2'=\varphi_{\tilde{w}_0}(\alpha_1+\alpha_2)=-\omega_2+\omega_3$ while 
$\mu_3'=\varphi_{\tilde{w}_0}(\alpha_1)$,
$\mu_1'=\varphi_{\tilde{w}_0}(\alpha_2)$ verify  
$\mu_3'+\mu_1'=-\omega_2+2 \omega_3 $. 

The set $\mathcal{A}$ of adapted convex orderings (commutation classes) is a 
small subset of $\mathcal{O}$ of all convex orderings, we provide
 the table below for small ranks :

\medskip

\hspace*{2cm}
\begin{tabular}{| c | c | c | c | c | c | c |}
\hline
 & $A_3$ & $A_4$ & $A_5$ & $A_6$ & $D_4$ & $D_5$ \\
\hline 
card $\mathcal{A}$ & $4$ & $8$ & $16$ & $32$ & $8$ & $16$ \\
\hline 
card $\mathcal{O}$ & $8$ & $62$ & $908$  & $24968$ & $182$ & $13198$ \\
\hline 
\end{tabular}
\bigskip

\textbf{Conjecture} $\tilde{w}_0$ is adapted to 
a quiver iff the mapping $\varphi_{\tilde{w}_0}$ 
is linear.

Computer testing shows 
linearity fails for all convex orderings in 
$\mathcal{O} \backslash \mathcal{A}$
for the Dynkin types in the table above. It is not clear to us
at the moment how linearity should be modified in order to 
describe $\varphi_{\tilde{w}_0}$ for arbitrary
commutation classes.

The set $\mathcal{C}_{\tilde{w}_0}$ may be stratified into levels 
the same way $\mathcal{R}_{\tilde{w}_0}$ is, namely the $i^{\mbox{\small th}}$
level $\mathcal{C}^i_{\tilde{w}_0}$ 
is the set of weights $\mu_k=s_{i_1} \ldots s_{i_k} (\omega_{i_k})$ 
with $i_k$ equal to $i$·

\medskip

\textbf{Lemma 2.6} Forall $i \in I$, the last element of $\mathcal{C}^i_{\tilde{w}_0}$
is $-\omega_{i^{*}}$.

\textbf{Proof :}  Consider the last appearance 
$s_{j_k}=s_i$ in $\tilde{w}_0$. All following reflections
$s_{j_m}$ with $m>k$ stabilize $\omega_i$. Hence
$s_{j_1} s_{j_2} \ldots s_{j_k} \omega_i=w_0 \omega_i=-\omega_{i^*}$ $\square$

\medskip

\textbf{Proposition 2.7} Let $\tilde{w}_0$ be adapted to $Q$, and $c$
the Coxeter element associated to $Q$. Then for all $i \in I$ : 
$$\mathcal{C}^i_{\tilde{w}_0}=
(c^{N_i-1} (-\omega_{i^*}) , \, c^{N_i-2} (-\omega_{i^*}), \ldots, 
\, c(-\omega_{i^*}), \, -\omega_{i^*}) $$ 
where $N_i=N_i(Q)$.

\textbf{Proof : } Recall $c=s_{j_1} s_{j_2} \ldots s_{j_n}$ with
$j_1$ a sink of $Q_1:=Q$, $j_2$ a sink of $Q_2:=S_1 Q$ and so on. As 
$(j_1, \, j_2, \, \ldots ,\,  j_n)$ is a permutation of 
$\{ 1, \ldots , \, n \}$, every arrow in $Q$ is inverted twice in
the sequence $S_{j_n} \ldots S_{j_1} Q$ so have $Q=Q_{n+1}$. 
Using Corollary 2.3 and $R_{n+1}=R_{1}=R$, we get 
$c \varphi_R=\varphi_R c$. By Theorem 2.4 we can  replace
$\varphi_R$ by $\varphi_{\tilde{w}_0}$ \linebreak[4]
so that : 
$$(*) \; \; c \, \varphi_{\tilde{w}_0}=\varphi_{\tilde{w}_0} \, c.$$

We already know by Auslander-Reiten theory the levels of 
$\mathcal{R}_{\tilde{w}_0}$ are given by the action of $c$. It 
remains to apply $(*)$ to get the statement of the proposition.
$\square$  

\textbf{Remark : } If $\tilde{w}_0$ is not adapted to a quiver $Q$, the levels
 $\mathcal{C}^i_{\tilde{w}_0}$ may not be given through 
the action of a Coxeter element. We do not have for the moment a Weyl
group model for the levels in general case.

\section{The Coxeter complex construction}

\hspace*{4mm} In this section we consider elements of 
$\mathcal{C}_{\tilde{w}_0}$ as one dimensional rays 
in the Coxeter complex $\Sigma$. The structure of $\Sigma$ allows 
 us to construct arrows, and associate a quiver $\Gamma_{\tilde{w}_0}$
with $\tilde{w}_0$. When $\tilde{w}_0$ is adapted to a quiver $Q$,
we show we recover $\Gamma_Q$.

Recall the geometry of the root system $\Phi$
is captured by the Coxeter complex $\Sigma$ (see \cite{brown} chapters 1-3, 
\cite{bourbaki} chapter V). 
Geometrically $\Sigma$ is the simplicial complex induced by the hyperplanes 
orthogonal to the positive roots. Maximal elements are 
the open chambers.  Recall a parabolic subgroup of $W$
is a subgroup $W^{J}=\langle s_j \rangle_{j \in J}$ with 
$J$ a subset of $I$. $\Sigma$ can be constructed out of $W$
as the complex of all left classes of all parabolic subgroups
partially ordered by inverse inclusion.

It would be more convenient for us to work in $\Sigma^{op}$,
the inversely ordered complex, as the partial order becomes
plain inclusion. Note $W^{(J)}:=\langle s_j \rangle_{j \in I \backslash J}$. 
The maximal elements (facets) of  $\Sigma^{op}$ are left classes of 
maximal parabolics $W^{(i)}$ (equivalently, one dimensional rays). 
Two facets $w W^{(i)}, \, w' W^{(j)}$ 
of $\Sigma^{op}$ are adjacent iff their intersection is 
non-empty, which in that case, it is a left class of $W^{(i,j)}$.

Let us recall now basic facts about parabolic subgroups, one may
consult for details \cite{humphreys} Chapter 1 or
\cite{geckpfeiffer} Chapters 1 and 2 (Geck and Pfeiffer  
work with right classes, but all the arguments can
be easily adapted to left classes). 

Fix $J \subset I$. $W^J$ is a itself a Coxeter group,
and as such is endowed with its proper length function
$l_J$.

\medskip

\textbf{Proposition 3.1} (\cite{geckpfeiffer} 1.2.10)
 Consider $w \in W^J$. Then $l_J(w)=l(w)$. Moreover,
if $\tilde{w}=s_{i_1} \ldots s_{i_m}$ is any reduced 
expression of $w$ inside $W$, then 
$\forall k, \; i_k$ is in $J$.

\medskip

\textbf{Proposition 3.2} (\cite{geckpfeiffer} 2.1.1)
Define 
$$ X_J=\{ w \in  W \mid l(ws)>l(w) \; \forall x \in J \}.$$

\textbf{a)} Every left class $w W^J$ contains a unique
element $u \in X_J$ and  $u$ is the unique element of minimal
length inside $w W^J$.

\textbf{b)} For any element $y \in W^J$, if we
write $y=u x$ with $x \in W^J$, then $l(y)=l(u)+l(x)$. 

\medskip

Consider $J \subset K \subset I$, and note by
$X_J^K$ the set of minimal representatives of
the left classes of $W^J$ in $W^K$

\medskip

\textbf{Lemma 3.3} ( \cite{geckpfeiffer} 2.2.5)
$$ X_J=X_K X_J^K.$$ 
 
Finally, $W^{(i)}$ is the stabilizer of $\omega_i$, thus 
a weight $\mu=w \omega_i$ may be identified with the left
class  $w W^{(i)}$. If $u$ is the minimal length element
inside $w W^{(i)}$ we shall say $u \omega_i$ is
the minimal presentation of $\mu$.

\medskip

\textbf{Lemma 3.4}  (\cite{BZschubert} 2.7) Let 
$\mu=w \omega_i$, and $\tilde{w}=s_{i_1} s_{i_2} \ldots s_{i_m}$
a reduced expression of $w$. 
The minimal representation of $\mu$ is obtained by erasing
from $\tilde{w}$ all factors $i_k$ for which
$\langle \omega_i, \,  \beta_k^{\vee} \rangle=0$, where  
$\beta_k=s_{i_m} s_{i_{m-1}} \ldots s_{i_{k+1}} (\alpha_{i_k})$.

Chamber weights can be seen as facets of $\Sigma^{op}$ by
identification above. We shall feel free to consider them
as appropriate either as weights or left classes.

\medskip

\textbf{Definition 3.5} 

Let $\tilde{w}_0$ be a reduced expression of $w_0$ in W a Weyl group of
finite type. The abstract Auslander-Reiten quiver 
$\Gamma_{\tilde{w}_0}$ associated with $\tilde{w}_0$ is defined as follows
\begin{itemize}
\item The set of vertices is the totally ordered set 
$\mathcal{C}_{\tilde{w}_0}=(\mu_1,  \, \ldots, \, \mu_N)$, 
seen as lying inside $\Sigma^{op}$.
\item There is an arrow  from  $\mu_j=w_j W^{(i_j)}$ to $\mu_k=w_k W^{(i_k)}$ iff
\begin{itemize}
\item[i)] $ \; \; i_j$ and $i_k$ are linked in Dynkin diagram;
\item[ii)]  $w_j W^{(i_j)}$ and $w_k W^{(i_k)}$ are adjacent in $\Sigma$;
\item[iii)] $j<k$.
\end{itemize}
\end{itemize}

We are forced to introduce $iii)$ as intersection 
of left classes is non-oriented.
Observe also that by considering only "Dynkin couples" in 
$i)$, $\Gamma_{\tilde{w}_0}$
actually depends only on the commutation class $\lbrack \tilde{w}_0 \rbrack$
of $\tilde{w}_0$. Definition 3.5 provides therefore 
a unique quiver $\tilde{\Gamma}_Q$ associated to the reduced
expressions adapted to a quiver $Q$. Our aim is to
show it is isomorphic to $\Gamma_Q$. We shall follow a combinatorial
approach, as suggested by the referee. 

\textbf{Example 3.1 :}
Consider the quiver 
$Q : \stackrel{1}{\bullet} \longrightarrow \stackrel{2}{\bullet} 
     \longleftarrow \stackrel{3}{\bullet}$
of type $A_3$. An adapted expression to $Q$ is 
$\tilde{w}_0=s_{i_1} s_{i_2} \ldots s_{i_6}$ with 
$\underline{i}=(2,\, 1, \, 3, \, 2 ,\, 3 ,\, 1)$.
$\Gamma_{\tilde{w}_0}$ is 

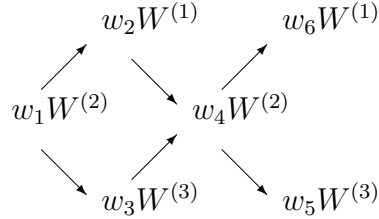
\begin{picture}(60,30)(-40,0)
\put(16,16){\vector(1,1){6}}
\put(28,20){\vector(1,-1){6}}
\put(40,16){\vector(1,1){6}}
\put(16,8){\vector(1,-1){6}}
\put(28,4){\vector(1,1){6}}
\put(40,8){\vector(1,-1){6}}
\put(12,12){$w_1 W^{(2)}$}
\put(24,24){$w_2 W^{(1)}$}
\put(24,0){$w_3 W^{(3)}$}
\put(36,12){$w_4 W^{(2)}$}
\put(48,24){$w_6 W^{(1)}$}
\put(48,0){$w_5 W^{(3)}$}
\end{picture}

\noindent where $w_k=s_{i_1} \ldots s_{i_k}, \; k=1 , \ldots 6$. 

Note for example
$$ \begin{array}{rl} 
w_1 W^{(2)} & =\{ s_2, \, s_2 s_1,\, s_2 s_3, \, s_2 s_1 s_3 \} \\
w_2 W^{(1)} & =\{ s_2 s_1, \, s_2 s_1 s_2 ,\, s_2 s_1 s_3, \, 
s_2 s_1 s_2 s_3, \, s_2 s_1 s_3 s_2, \, s_2 s_1 s_2 s_3 s_2 \}
\end{array}$$ 
are adjacent since 
$w_1 W^{(2)} \cap w_2 W^{(1)}=\{ s_2 s_1, \, s_2 s_1 s_3 \}
(=s_2 s_1 W^{(1,2)})$, while $w_1 W^{(2)}$ and $w_6 W^{(1)}$ 
are not since $w_1 W^{(2)} \cap w_6 W^{(1)}=\emptyset$.  We
leave it to the reader to verify that there is an arrow 
between vertices on adjacent levels (=Dynkin couples) exactly
when intersection of left classes is non empty.

Let $\tilde{w}=s_{i_1} s_{i_2} \ldots s_{i_m}$ a
reduced expression of an arbitrary element 
$w \in W,$ \linebreak[4] 
$w \neq e$. We can define a partial
set of chamber weights 
$\mathcal{C}_{\tilde{w}}=(\mu_1, \, \ldots , \, \mu_n)$
with \linebreak[4]
$\mu_j=s_{i_1} s_{i_2} \ldots s_{i_{j-1}} (\mu_{i_j})$, 
and a partial quiver $\Gamma_{\tilde{w}}$ as above.

\medskip

\textbf{Proposition 3.6}
Assume $\tilde{w}$ is alternating.
Then combinatorial the description of Proposition $1.2$ is valid
for $\Gamma_{\tilde{w}}$. 
There is an arrow $\mu_j \longrightarrow \mu_k$ 
if and only if 
\begin{itemize}
\item[i)] $i_j$ and $i_k$ are linked in the Dynkin diagram $D$;
\item[ii)] $j<k$;
\item[iii)] For $j < p < k$, $i_p \neq i_j$.
\end{itemize}

\textbf{Proof : } 
We proceed by induction on $l(w)$. If $l(w)=1$ there is nothing
to prove. Suppose the proposition is valid for an 
alternating reduced expression $\tilde{w}=s_{i_1} \ldots s_{i_m}$, and
consider an alternating $\tilde{w}'=s_r \tilde{w}$ with 
$l(w')=l(w)+1$. Note 
$\mathcal{C}_{\tilde{w}'}=(\mu_1', \, \ldots, \mu'_{m+1})$.  

For any $x, \, y \in W, \;  x W^{(i)} \cap y W^{(j)} \neq \emptyset 
\Longleftrightarrow s_r x W^{(i)} \cap s_r y W^{(j)} \neq \emptyset$. 
Using induction hypothesis for $\tilde{w},$
the combinatorial description is valid for arrows 
$\mu_j' \longrightarrow \mu_k'$  
inside $\Gamma_{\tilde{w}'}$ with $j \geq 2$. 
It remains to study arrows out of $\mu_1'$, that is 
intersections
$s_r W^{(r)} \cap s_r s_{i_1} \ldots s_{i_k} W^{(i_k)}$. By
the above, this amounts to consider intersections 
$W^{(r)} \cap s_{i_1} \ldots s_{i_k} W^{(i_k)}$. An arrow
$\mu_1' \longrightarrow \mu_{k+1}'$ exists in $\Gamma_{\tilde{w}'}$
iff such an intersection is non-empty.
We shall consider  $ 1 \leq k \leq m$ with $i_k$ linked to $r$ in $D$
and put $w_k=s_{i_1} \ldots s_{i_k}$.

Let us show first all arrows $\mu_1' \longrightarrow \mu_{k+1}'$ 
given by the criterion appear in $\Gamma_{\tilde{w}'}$.  
If $k$ verifies for all $1<p<k, \; i_p \neq r$, as also $i_k \neq r$, 
we have $w_k \in W^{(r)}$. 
Hence  $W^{(r)} \cap w_k W^{(i_k)} \neq \emptyset$ as expected.

Conversely, let us show there are no more arrows 
$\mu_1' \longrightarrow \mu_{k+1}'$ in $\Gamma_{\tilde{w}'}$
than those given by the criterion.
We suppose therefore there is $p$, $1<p<k$ such that $i_p=r$. 
We can take $p$ maximal with that property. $\tilde{w}$ is
alternating by hypothesis, hence 
$i_{p+1}, i_{p+2} \ldots i_{k-1}$ are different from 
$r$ and $i_k$. By 
definition of the simple reflections we have the following
coefficients ( $\alpha_{i_p}=\alpha_r)$) :
$$ (*)  
\left\{
\begin{array}{rl}
\lbrack s_{i_{k-1}} s_{i_{k-2}} \ldots s_{i_{p+1}} (\alpha_{i_p}) \, : \, 
\alpha_r \rbrack& =1 \\
\lbrack s_{i_{k-1}} s_{i_{k-2}} \ldots s_{i_{p+1}} (\alpha_{i_p}) \, : \, 
\alpha_{i_k} \rbrack& =0 
\end{array}
\right. $$
Now $s_{i_k} (\alpha_{i_p})=\alpha_r+\alpha_{i_k}$, so we obtain from
(*) :
$$ \lbrack s_{i_{k}} s_{i_{k-1}} \ldots s_{i_{p+1}} (\alpha_{i_p}) \, : \, 
\alpha_{i_k} \rbrack > 0$$
and therefore 
$( \omega_{i_k}, \, 
s_{i_{k}} s_{i_{k-1}} \ldots s_{i_{p+1}} (\alpha_{i_p}^{\vee}) ) \neq 0$.  
By Lemma 3.4, $s_{i_p}=s_r$ appears in
 the reduced expression $\tilde{u}$ for minimal representative $u$ 
of $w_k W^{i_k}$ obtained by erasing procedure out of  $\tilde{w}_k$.

Consider a left class of $W^{(r,i_k)}$ inside $w_{i_k} W^{(i_k)}$.
By Lemma 3.3, its minimal representative is of the form
$y=u x$ with $x$ a minimal representative of a class of  $W^{(r,i_k)}$
inside    $W^{(i_k)}$. Take any reduced expression $\tilde{x}$ of $x$.
$\tilde{y}=\tilde{u} \tilde{x}$ is a reduced expression of $y$, in
which $s_r$ appears. By Proposition 3.1 this means $y \notin W^{(r)}$,
hence $y W^{(r,i_k)}$ is not included in $W^{(r)}$. We see
no left class of $W^{(r,i_k)}$ can be contained in both
$W^{(r)}$ and $w_k W^{(i_k)}$, so there is no
arrow in $\Gamma_{\tilde{w}'}$ from $\mu_1'$ towards $\mu_{i_{k+1}}'$ 
$\square$
 
\textbf{Remark :  } 
 Experimental data suggests the statement of 
the proposition is actually true in general, without the
assumption that $\tilde{w}$ is alternating (see example 3.2 below).
For this reason we insist that the third condition in the combinatorial
characterization be $i_p \neq i_j$
while $i_p \neq i_j, \, i_k$ would have been sufficient 
for the adapted case.

\medskip

\textbf{Theorem 3.7} (Isomorphism)

Let $Q$ be a quiver of type $ADE$, and $R$ its Ringel form.
Then $\varphi_R$ induces isomorphism of oriented graphs 
of $\Gamma_Q$ onto
$\tilde{\Gamma}_Q$.

\textbf{Proof :}
Let $\tilde{w}_0$ be a reduced expression adapted to $Q$
Combining Theorems 1.1 and 2.4, we get a  correspondence 
between vertices of $\Gamma_Q$ and $\tilde{\Gamma}_Q$
given by $\varphi_R$. By Lemma 1.4, $\tilde{w}_0$
is alternating, hence applying  Propositions 1.2 and 3.6, 
$\varphi_R$ extends to an isomorphism of oriented graphs.

Theorem 3.7 justifies the use of the term Abstract Auslander-Reiten 
quiver for $\Gamma_{\tilde{w}_0}$. There is a finite
dimensional algebra, behind the commutation class
corresponding to $Q$, namely $\Cbold Q$.  Definition 3.5 arises the question
of associating a finite dimensional algebra with 
any commutation class. One  aim would be 
controlling the corresponding parameterizations of
the canonical basis (see \cite{reineke} for some adapted cases). 
Another application of finite dimensional algebra techniques
in general reduced expressions setting,   
might  be solving the problem of realizability of 
the set of  commutation classes with braid moves (see \cite{bedard}) 
as a $1$-skeleton of a polytopal complex.
This problem is also known as realizability of
second Bruhat order \cite{ziegler}, we point for instance to 
\cite{felsnerziegler} as representing the current trend 
of research concerning this problem. 
We observe for the moment, even the number of commutation classes in
$A_n$ case is unknown. 

\textbf{ Example 3.2 : } 
We have seen linearity fails for non-adapted cases. The same applies
for the ``mesh'' structure coming from $\Zbold D$. 
Take $\tilde{w}_0=s_2 s_1 s_2 s_3 s_2 s_1$. One may check
$s_2 s_1 W^{(1)} \cap s_2 s_1 s_2 s_3 s_2 W^{(2)}=s_2 s_1 s_3 s_2 W^{(1,2)}$.
$\Gamma_{\tilde{w}_0}$ is given by 

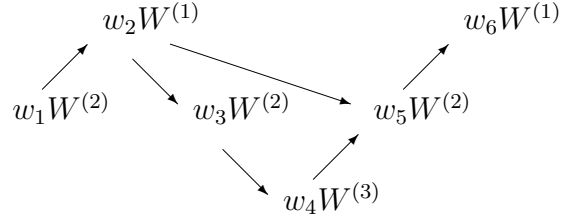
\vspace*{-2.0cm} \hspace*{2cm}
\begin{picture}(80,50)
\put(16,16){\vector(1,1){6}}
\put(28,20){\vector(1,-1){6}}
\put(40,8){\vector(1,-1){6}}
\put(52,4){\vector(1,1){6}}
\put(64,16){\vector(1,1){6}}
\put(33,22){\vector(3,-1){24}}
\put(12,12){$w_1 W^{(2)}$}
\put(24,24){$w_2 W^{(1)}$}
\put(36,12){$w_3 W^{(2)}$}
\put(48,0){$w_4 W^{(3)}$}
\put(60,12){$w_5 W^{(2)}$}
\put(72,24){$w_6 W^{(1)}$}
\end{picture}
\bigskip
\hfill \break
We see there are two arrows from $s_2 s_1 W^{(1)}$
towards the second level. Observe also $\tilde{w}_0$ is no longer 
alternating, erasing $s_3$ we get $s_2 s_1 s_2 s_ 2 s_1$. 
The combinatorial description 
of Proposition 3.6 still applies.

\section{The $A_n$ case : wiring diagrams}

\hspace*{5mm}  We describe in this section $A_n$ case.
We use Young tableaus (columns) to describe 
$\tilde{\Gamma}_Q$. On the 
theoretical level, we describe briefly the correspondence
of abstract Auslander-Reiten quivers  with wiring diagrams
used in \cite{BFZ}.

Linearity allows us to compute easily the vertices of $\Gamma_{\tilde{w}_0}$,
however, we are unable for the moment to provide an generalized 
adjacency criteria inside $\Sigma^{op}$, which would give in 
a straight-forward manner the arrows. There are in the literature 
case by case combinatorial descriptions of $\Sigma$ according to 
Dynkin type, we restrict ourselves
here to  the simplest case : $A_n$.  

All fundamental representations $E(\omega_i)$ of a simple complex lie algebra
of type $A_n$ are minuscule, so that $W \omega_i$ is the set of weights
of $E(\omega_i)$. Furthermore as $E(\omega_i)$ is isomorphic to  
$\bigwedge^i E(\omega_1)$, it has a basis of weight vectors indexed 
by all Young columns of size $i$, filled in a strictly increasing manner
with indices in $\{ 1 , \ldots , n+1 \}$. We get therefore an identification
of  $W \omega_i$ (equivalently the set of left classes of $W^{(i)}$) 
with multi-indices $J=(j_1, j_2 \ldots j_i), \;  
1 \leq j_1 <  j_2 \ldots < j_i \leq n$.
Conversely, given a strictly increasing multi-index 
$J$ of cardinality $i$ (that is
a Young column of size $i$), one may recover its weight as the sum of weights
of its boxes 
$$wt(J)=\underset{k=1}{\overset{i}{\sum}} wt( \fbox{$j_k$}),$$
with $wt( \fbox{$j_k$})=\omega_k-\omega_{k-1}$, as such a box
 corresponds to a weight vector inside $E(\omega_1)$
(by convention $\omega_0=\omega_{n+1}=0$).

In $A_n$ case, $W$ is the symmetric group $\mathcal{S}_{n+1}$, so that simple
reflections are the transpositions $(i \, i+1)$, and action on weight
orbits becomes usual $\mathcal{S}_{n+1}$ action on sets of indices.
Let us note for two integers $a<b$ the interval 
$(a,\, a+1, \, \ldots, \, b)$ by $\lbrack a , \, b \rbrack$.

\medskip

\textbf{Proposition 4.1} Fix $i$ between $1$ and $n$, and consider 
two left cosets $w W^{(i)}, \; w' W^{(i+1)}$, with $J, \; J'$ 
corresponding multi-indices. Then $w W^{(i)}$ and $w' W^{(i+1)}$ are adjacent
in $\Sigma^{op}$ if and only if $ J \subset J'$.

\textbf{Proof : }  The closure of the dominant chamber
is a fundamental domain for the action of $W$. Furthermore
pairs $(W^{(i)}, \, w' W^{(i+1)})$ correspond to left 
classes of $W^{(i,i+1)}$ inside $W^{(i)}$ which form
a single $W^{(i)}$ orbit. Hence, the adjacent pairs 
$(w W^{(i)}, \, w' W^{(i+1)})$ in $\Sigma^{op}$ form one
orbit under diagonal action, that of $( W^{(i)},W^{(i+1)})$.

The multi-indices corresponding to $W^{(i)}$ and $W^{(i+1)}$ 
are $\lbrack 1, \, i \rbrack$ and $\lbrack 1 , \, i+1 \rbrack$. 
As Weyl group action preserves inclusion of multi-indices we get
from the discussion above that $w W^{(i)}, \, w' W^{(i+1)}$
adjacent implies $J \subset J'$.

For the converse we need to show any couple 
$(J, \, J')$ with $J \subset J'$ can be obtained by 
diagonal action from 
$(\lbrack 1, \, i \rbrack, \, \lbrack 1, \, i+1 \rbrack)$.
Since all multi-indices of cardinality $i$ form a single Weyl group
orbit, it is enough to consider couples of the form 
$(\lbrack 1, \, i \rbrack , \, J')$. Since 
$\lbrack 1 , \, i \rbrack \subset J'$, 
we must have $J'=\lbrack 1 , \, i \rbrack \cup \{ k \}$ with
$k \geq i+1$. If $k>i+1$, 
then $s_{i+1} s_{i+2} \ldots s_{k-1} J'=\lbrack 1, i+1 \rbrack$
while $s_{i+1} s_{i+2} \ldots s_{k-1}$ stabilizes 
 $\lbrack 1, \, i \rbrack$. Thus 
$$( \lbrack 1 , \, i \rbrack, J')=s_{k-1} s_{k-2} \ldots s_{i+1} 
( \lbrack 1 , \, i \rbrack, \; \lbrack 1 , \, i+1 \rbrack)
\; \; \; \; \; \; \; \; \square$$

Let us now describe $\tilde{\Gamma}_{Q}$ in terms of column tableaus. 
It seems somewhat simpler to use the classical method 
of Coxeter translation, instead of linearity.

As $W \cong \mathcal{S}_{n+1}$, the Coxeter
element $c$ defined by $Q$ is a $n+1$ cycle. Let us fix a cycle writing 
$\lbrack c \rbrack=(a_1 \, a_2 \, \ldots \, a_{n+1})$ by
imposing $a_{n+1}=n+1$. In order to work inside $\lbrack c \rbrack$, we shall use
$\overline{m}$  for $m \; \mbox{mod} \; n+1$. Finally let 
$\lbrack c \rbrack_{\overline{m},i}$ be the multi-index obtained by
ordering the "interval" 
$\{ a_{\overline{m}},\, a_{\overline{m+1}}, \, \ldots, \, a_{\overline{m+i-1}} \}$.

\medskip

\textbf{Lemma 4.2 } Define for $i \in I, \; \overline{m_i}$ 
recursively by $\overline{m}_1=n+1$ and for $i \geq 2$
$$ \overline{m}_i= \left\{ 
\begin{array}{ll}
\overline{m}_{i-1} & \mbox{ if } 
\stackrel{i^*}{\bullet} \longrightarrow \stackrel{i^*+1}{\bullet}
\mbox{ in } Q \\
\overline{m}_{i-1}-1 & \mbox{ if } 
\stackrel{i^*}{\bullet} \longleftarrow \stackrel{i^*+1}{\bullet}
\mbox{ in } Q
\end{array} \right.$$  
Then $[c]_{\overline{m}_i,i}$ is the multi-index 
$(i^*+1, \, i^*+2, \, \ldots , \, n+1)$.

\textbf{Proof : } We proceed by induction on $n$. For case $n=2$ there is
nothing to prove. Consider $Q$ of type $A_n$, and let $Q'$ be its
sub-quiver with vertices $1$ to $n-1$, with $c'$ the corresponding
Coxeter element, and $m'$ defined for $c'$ as above. 

\underline{Case 1} : Q is obtained from $Q'$ by adding 
$\stackrel{n-1}{\bullet} \longrightarrow \stackrel{n}{\bullet}$. $n$ is a 
sink, hence  $c= s_n c=(n \, n+1) c'$. If
$c'=(b \ldots a \, n )$ then $c=(b \ldots a n+1 n)$.
We have therefore $\overline{m}_2=\overline{m}_1$.

\underline{Case 2} : Q is obtained from $Q'$ by adding 
$\stackrel{n-1}{\bullet} \longleftarrow \stackrel{n}{\bullet}$. $n$ is
a source, hence $c=c' s_n=c' (n \, n+1)$. If 
$c'=(b \ldots a \, n )$ then $c=(b \ldots a n n+1)$.
We get $\overline{m}_2=\overline{m_1-1}$.

In both cases, the cycle writing of $c$ is obtained 
from that of $c'$ by replacing $n$ by either $n \, n+1$
or $ n+1 \, n$. The validity of the Lemma 
for $c'$ in cases $i=2 \ldots n-1$ implies its validity
for $c$ in cases $i=3 \ldots n$. Also, the recursive 
formula for
$\overline{m}_i$ for $i \leq 3$ is the same as that for
$\overline{m}'_{i-1}$ $\square$

We see $c$ has a ``segment property'' : for any $j$, 
indices $j, j+1, \ldots n+1$ form a ``connected'' band  (modulo $n+1$)
in the cycle writing of $c$. This property leads to the
description of levels of $\tilde{\Gamma}_Q$.
Let use note by $J^{(i)}_j(Q)$ the
multi-index of the $j^{\mbox{{\small th}}}$ vertex from the right 
(starting with $0$) 
on the $i^{\mbox{{\small th}}}$ level of $\tilde{\Gamma}_Q$. 

\medskip

\textbf{Proposition 4.3}
Consider $Q$ and fix $i \in I$, then 

\textbf{a)} (\cite{bedard} Corollary 2.20 (a)) 
The cardinality $N_i(Q)$ of the $i^{\mbox{\small th}}$ level of 
$\tilde{\Gamma}_{Q}$ is \linebreak[4]
$(n+1+a_i-b_i)/2$ where $a_i$ (respectively $b_i$) is the
number of arrows in the unoriented path from $i$ to $i^*$ 
that are directed towards $i$ (respectively $i^*$).

\textbf{b} For $j=0,\, 1, \, \ldots, \, N_i(Q)-1$ one has
$$ J_{j}^{(i)}(Q)=\lbrack c \rbrack_{\overline{m_i-j},i}$$

\textbf{Proof of b) : } Fix $i \in I$. By Lemma 2.6, 
the first weight from the right
on the $i^{\mbox{\small th}}$ level of $\tilde{\Gamma}_Q$ is
$-\omega_{i^*}$. It corresponds to the multi-index 
$(i^*+1, \, i^*+2, \ldots \, n+1)$. By the Lemma above, 
the proposition is true for $j=0$.

For $ 1 \leq j \leq N_i(Q)-1$, we use Proposition 2.7 :
$$ \begin{array}{rl}
J_j^{(i)}(Q) &= c^{-j} J_{0}^{i}(Q) \\
             &= c^{-j} \lbrack c \rbrack \\
             &= \lbrack c \rbrack_{\overline{m_i-j},i} \square
\end{array} $$

We see getting $\tilde{\Gamma}_Q$ out of $Q$ requires two
operations : computing the cycle writing of the Coxeter 
element associated to $Q$, and the cardinalities $N_i(Q)$.
Both involve working directly with $Q$. Once vertices are 
known, getting the arrows is immediate due to Lemma 4.1.

\textbf{Example 4.1 : }
Consider the quiver $Q : \; \stackrel{1}{\bullet} \longrightarrow 
\stackrel{2}{\bullet} \longleftarrow \stackrel{3}{\bullet}$ of type $A_3$. 
$1^*=3$ so that $a_1=1, \, b_1=1$ 
hence $N_1(Q)=2$. We obtain in a similar manner
$N_2(Q)=2, \; N_3(Q)=2$. The Coxeter element $c$ is the cycle
$(2 1 3 4)$. The levels of $\tilde{\Gamma}_Q$ are given respectively
by $\tilde{\Gamma}^1_{Q} \; : \; (\cdot \cdot 3 \cdot), (\cdot \cdot \cdot 4);
\; \tilde{\Gamma}^2_Q \; : \; (\cdot 1 3 \cdot), \; (\cdot \cdot 3 4); \; 
\tilde{\Gamma}^3_Q \; : \; (\cdot 1 3 4), \; (2 \cdot 3 4)$. 
Considering inclusion between sets of indices we get the arrows : 

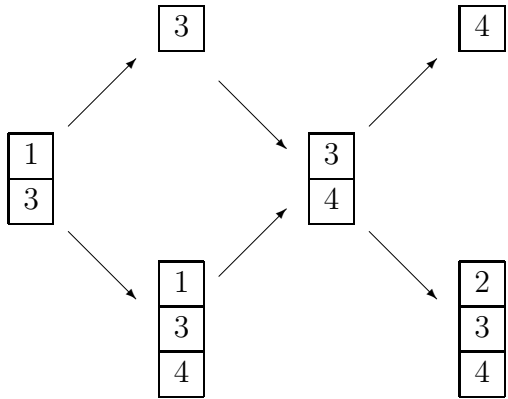
\begin{picture}(80,50)
\put(27,27){\vector(1,1){9}}
\put(47,33){\vector(1,-1){9}}
\put(67,27){\vector(1,1){9}}
\put(27,13){\vector(1,-1){9}}
\put(47,7){\vector(1,1){9}}
\put(67,13){\vector(1,-1){9}}
\put(19,14){\line(1,0){6} }
\put(19,14){\line(0,1){6} }
\put(19,20){\line(1,0){6} }
\put(25,14){\line(0,1){6} }
\put(21,16){3}
\put(19,20){\line(1,0){6} }
\put(19,20){\line(0,1){6} }
\put(19,26){\line(1,0){6} }
\put(25,20){\line(0,1){6} }
\put(21,22){1}
\put(39,37){\line(1,0){6} }
\put(39,37){\line(0,1){6} }
\put(39,43){\line(1,0){6} }
\put(45,37){\line(0,1){6} }
\put(41,39){3}
\put(39,-9){\line(1,0){6} }
\put(39,-9){\line(0,1){6} }
\put(39,-3){\line(1,0){6} }
\put(45,-9){\line(0,1){6} }
\put(41,-7){4}
\put(39,-3){\line(1,0){6} }
\put(39,-3){\line(0,1){6} }
\put(39,3){\line(1,0){6} }
\put(45,-3){\line(0,1){6} }
\put(41,-1){3}
\put(39,3){\line(1,0){6} }
\put(39,3){\line(0,1){6} }
\put(39,9){\line(1,0){6} }
\put(45,3){\line(0,1){6} }
\put(41,5){1}
\put(59,14){\line(1,0){6} }
\put(59,14){\line(0,1){6} }
\put(59,20){\line(1,0){6} }
\put(65,14){\line(0,1){6} }
\put(61,16){4}
\put(59,20){\line(1,0){6} }
\put(59,20){\line(0,1){6} }
\put(59,26){\line(1,0){6} }
\put(65,20){\line(0,1){6} }
\put(61,22){3}
\put(79,-9){\line(1,0){6} }
\put(79,-9){\line(0,1){6} }
\put(79,-3){\line(1,0){6} }
\put(85,-9){\line(0,1){6} }
\put(81,-7){4}
\put(79,-3){\line(1,0){6} }
\put(79,-3){\line(0,1){6} }
\put(79,3){\line(1,0){6} }
\put(85,-3){\line(0,1){6} }
\put(81,-1){3}
\put(79,3){\line(1,0){6} }
\put(79,3){\line(0,1){6} }
\put(79,9){\line(1,0){6} }
\put(85,3){\line(0,1){6} }
\put(81,5){2}
\put(79,37){\line(1,0){6} }
\put(79,37){\line(0,1){6} }
\put(79,43){\line(1,0){6} }
\put(85,37){\line(0,1){6} }
\put(81,39){4}
\end{picture}

\hspace*{1cm}
\bigskip
\bigskip

Reduced expressions $\tilde{w}_0$ of the longest word in $W$ play
a central role in the theory of the canonical basis
\cite{kashbase} \cite{lustbase} of the positive part
$U_q(\nbold^+)$ of the quantized enveloping algebra
$U_q(\gbold)$. As this base is difficult to compute, 
Lusztig's model consists in approximating it by a basis
composed of PBW monomials. Every convex ordering of $\Phi^+$
(or commutation class) gives rise to a particular PBW base,
which provides an indexation of the canonical basis. The
control over these parameterizations involves a geometrical
interpretation of $\tilde{w}_0$, the wiring diagram 
$\mathcal{WD}(\tilde{w}_0)$.  We refer to \cite{BFZ} \S 2.3 for the 
\linebreak[4] details, 
and reduce the presentation here to example 4.2.
Our conventions are ``opposite'' (twisted by $*$) 
compared to those of \cite{BFZ}.

Given $\tilde{w}_0=s_{i_1} s_{i_2} \ldots s_{i_N}$ in 
$W \cong \mathcal{S}_{n+1}$, 
the wiring diagram $\mathcal{WD}(\tilde{w}_0)$ is
an arrangement of $n+1$ pseudo-lines of the plane $\Rbold^2$, 
going from $-\infty$ to $+\infty$.
The pseudo-lines cross pairwise in points $x_1, \ldots x_N$, situated
in levels (top to bottom) in accordance
with $\tilde{w}_0$, namely the crossing $x_k$ 
is on the level $i_k$. We get
a braid reversing the order of the lines, 
and which partitions $\Rbold^2$ into zones
(or chambers). One numbers the pseudo-lines from $1$ to $n+1$
as in the example (at $x$ coordinate $-\infty$ : from top to
bottom). This gives a labelling of each zone $Z$ by
the set $J(Z)$ of indices of lines passing above it.
Each zone has a distinct label. 
Observe \cite{BFZ} geometrically $\mathcal{WD}(\tilde{w}_0)$ depends
only on the commutation class of $\tilde{w}_0$, and there
is a $1 : 1$ bijection between the set of commutation
classes of $A_n$ and possible wiring diagrams with $n+1$ strands.   
\medskip  

\textbf{Example 4.2 : } Take $\tilde{w}_0=s_2 s_1 s_3 s_2 s_3 s_1$
 (type $A_3$). $\mathcal{WD}(\tilde{w}_0)$ is : 

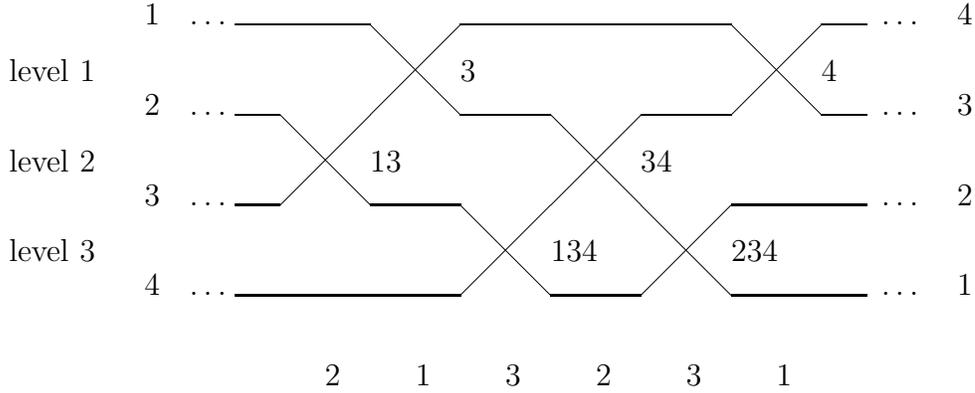
\vspace*{-4cm}
\hspace*{3cm}
\begin{picture}(120,96)
\put(-12,12){$4$}
\put(-12,24){$3$}
\put(-12,36){$2$}
\put(-12,48){$1$}
\put(0,12){\line(1,0){6}}
\put(0,24){\line(1,0){6}}
\put(0,36){\line(1,0){6}}
\put(0,48){\line(1,0){6}}
\put(-6,12){$\ldots$}
\put(-6,24){$\ldots$}
\put(-6,36){$\ldots$}
\put(-6,48){$\ldots$}
\put(96,12){$1$}
\put(96,24){$2$}
\put(96,36){$3$}
\put(96,48){$4$}
\put(78,12){\line(1,0){6}}
\put(78,24){\line(1,0){6}}
\put(78,36){\line(1,0){6}}
\put(78,48){\line(1,0){6}}
\put(86,12){$\ldots$}
\put(86,24){$\ldots$}
\put(86,36){$\ldots$}
\put(86,48){$\ldots$}
\put(-30.,16.5){level $3$}
\put(-30.,28.5){level $2$}
\put(-30.,40.5){level $1$}
\put(12,0){$2$}
\put(24,0){$1$}
\put(36,0){$3$}
\put(48,0){$2$}
\put(60,0){$3$}
\put(72,0){$1$}
\put(6,36){\line(1,-1){12} }
\put(6,24){\line(1,1){12} }
\put(6,12){\line(1,0){12} }
\put(6,48){\line(1,0){12} }
\put(18,48){\line(1,-1){12} }
\put(18,36){\line(1,1){12} }
\put(18,12){\line(1,0){12} }
\put(18,24){\line(1,0){12} }
\put(30,24){\line(1,-1){12} }
\put(30,12){\line(1,1){12} }
\put(30,36){\line(1,0){12} }
\put(30,48){\line(1,0){12} }
\put(42,36){\line(1,-1){12} }
\put(42,24){\line(1,1){12} }
\put(42,12){\line(1,0){12} }
\put(42,48){\line(1,0){12} }
\put(54,24){\line(1,-1){12} }
\put(54,12){\line(1,1){12} }
\put(54,36){\line(1,0){12} }
\put(54,48){\line(1,0){12} }
\put(66,48){\line(1,-1){12} }
\put(66,36){\line(1,1){12} }
\put(66,12){\line(1,0){12} }
\put(66,24){\line(1,0){12} }
\put(18,28.5){$13$}
\put(30,40.5){$3$}
\put(42,16.5){$134$}
\put(54,28.5){$34$}
\put(66,16.5){$234$}
\put(78,40.5){$4$}
\end{picture}
\bigskip

Consider the set of zones bounded to the left by a 
crossing point. These are in bijection with crossing points $x_k$.
Let us note therefore by $Z_k$ the zone whose leftmost point
is $x_k$ in $\mathcal{WD}(\tilde{w}_0)$, and associate to
it the level of $x_k$. 

\medskip

\textbf{Lemma 4.4}  

\textbf{a)} The label  $J(Z_k)$ is $w_k \cdot ( 1,2,\ldots, \, i_k)$
where  $w_k=s_{i_1} s_{i_2} \ldots s_{i_k}$.

\textbf{b)} Two chambers $Z_j$ on the $i^{\mbox{\small th}}$ level
and $Z_k$ on the $(i+1)^{\mbox{\small th}}$ level are adjacent 
if and only if a pseudo-line $L_m$ passes between them, that is
if and only if $J(Z_k)=J(Z_j) \cup \{ m \}$.

\textbf{Proof : } Pass a vertical line $p_k$ through each crossing $x_k$.
This divides the plane into vertical bands 
$\mathcal{V}_0, \, \mathcal{V}_1, \ldots ,\, \mathcal{V}_N$. 
$\mathcal{V}_0$ lies between $-\infty$ and $p_1$, 
for $k=1, \ldots , N-1$, $\mathcal{V}_k$ lies between 
$p_{k}$ and $p_{k+1}$, and $\mathcal{V}_N$ lies between 
 $p_N$ and $+\infty$.

Note by $J_k$ the ordered set of indices of pseudo-lines
passing through $\mathcal{V}_k$, read top to bottom. Initially
$J_0=(1, \, 2,\ldots,\, n)$. If we note $i=i_{k}$, the crossing
$x_{k}$ \linebreak[4] consists 
in interchanging the lines in positions $i$ and $i+1$.
Suppose \linebreak[4]
$ J_{k-1}=(w_{k-1}(1), \, w_{k-1}(2) , \, \ldots, 
w_{k-1}(i), \, w_{k-1}(i+1), \ldots, \, w_{k-1}(n))$, then
$$ \begin{array}{rl}
J_{k} &=(w_{k-1}(1), \, w_{k-1}(2) , \, \ldots, 
w_{k-1}(i+1), \, w_{k-1}(i), \ldots, \, w_{k-1} (n))  \\
&= (w_{k-1} s_i (1), \, w_{k-1} s_i (2) , \, \ldots, 
w_{k-1} s_i (i), \, w_{k-1}s_i (i+1), \ldots, \, w_k s_i (n)) \\
&=(w_k(1), \, w_k(2) \ldots w_k(n)). 
\end{array}$$
Proceeding by induction we get for 
$k=1 \ldots N, \; J_k=(w_k(1), \, w_k(2), \ldots, \, w_k(n))$.
Part  a) follows then by observing for a zone $Z_k$, the label $J(Z_k)$ 
consists in the first $i_k$ indices of $J_k$. 

Part b) follows directly from our definition of labelling of 
the zones in $\mathcal{WD}(\tilde{w}_0)$ $\square$

We may consider labels of zones $Z_k$ as Young columns. Lemma
4.4 a) shows this establishes a bijection with $\mathcal{C}_{\tilde{w}_0}$
(thus the term chamber weights used in \cite{BZschubert}).
Comparing Lemma 4.1 with Lemma 4.4 b), we see arrows in $\Gamma_{\tilde{w}_0}$
correspond to adjacencies of zones in $\mathcal{WD}(\tilde{w}_0)$.
We invite the reader to compare Examples 4.1 and 4.2.
Observe the non-left bounded zones may be considered as trivial, as
they appear in any wiring diagram, and always have the same label.
We therefore have a natural $1 : 1$ correspondence 
between the set of abstract Auslander-Reiten quivers and
that of wiring diagrams. 

\textbf{Remark : } Wiring diagrams were generalized to
$B_n, \,C_n$ cases in \cite{BZschubert}. No such construction
exists in $D_n$ or $E_n$ cases.

In $A_n$ case, we may summarize our results in the following figure : 

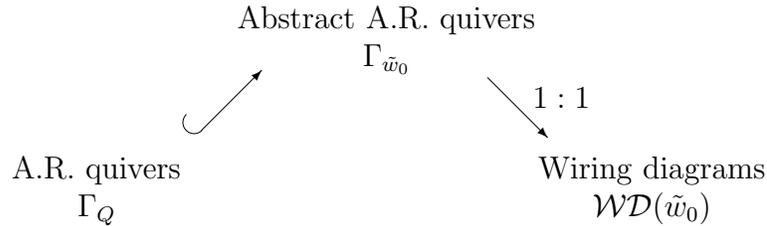
\vspace*{-1cm} \hspace*{2.3cm}
\begin{picture}(120,40)
\put(0,0){$\begin{array}{c} \mbox{A.R. quivers} \\ \Gamma_Q \end{array}$}
\put(27,9){\vector(1,1){8}}
\put(26,10){\arc{3}{7.06}{3.93}}
\put(30,20){$\begin{array}{c} \mbox{Abstract A.R. quivers} \\ \Gamma_{\tilde{w}_0} \end{array}$}
\put(65,16){\vector(1,-1){8}}
\put(71,12){$1:1$}
\put(70,0){$\begin{array}{c} \mbox{Wiring diagrams} \\ \mathcal{WD}(\tilde{w}_0) \end{array}$}   
\end{picture}
\bigskip
\bigskip

We obtain another instance of the link between theories  
of finite dimensional algebras and 
of the canonical basis of $U_q(\gbold)$.

\begin{center}
{ \small S. Zelikson, Laboratoire S.D.A.D.,  unit\'{e} CNRS  FRE 2272} \\ 
{ \small D\'{e}partement de Math\'{e}matiques} \\ 
{ \small Universit\'{e} de Caen B.P. 5186 } \\
{ \small 14032 Caen Cedex FRANCE} \\
{ \small E-mail : Shmuel.Zelikson@math.unicaen.fr}
\end{center}

\end{document}